\documentclass[11pt, a4paper, twoside]{amsart}
\usepackage{mathrsfs, tabularx}

\usepackage[all]{xypic}
\usepackage {amsthm, amsmath, amssymb, amscd, float, latexsym, times, epic, eepic}

\theoremstyle{plain}
\newtheorem{theorem}{Theorem}[section]
\newtheorem{lemma}[theorem]{Lemma}
\newtheorem{proposition}[theorem]{Proposition}
\theoremstyle{definition}
\newtheorem{definition}[theorem]{Definition}

\newtheorem{remark}[theorem]{Remark}

\newcommand {\Set}[1] {\mathbb{#1}}
\newcommand{\setR}[0]{\Set{R}}
\newcommand{\setC}[0]{\Set{C}}

\newcommand{\bigMap}[0]{F}

\newcommand {\proofBox}[0]{\hfill $\Box$ }
\newcommand {\commentOut}[1]{}
\newcommand {\wt}[0]{\widetilde}
\newcommand {\emi}[0]{\emph{(i)} }
\newcommand {\emii}[0]{\emph{(ii)} }
\newcommand {\emiii}[0]{\emph{(iii)} }

\newcommand {\proofread}[1]{}

\newcommand{\pd}[2]{\frac{\partial #1}{\partial #2}}

\newcommand{\pdd}[3]{\frac{\partial^2 #1}{\partial #2\,\partial #3}}

\newcommand{\vfield}[1]{{\mathfrak X}( #1)}
\newcommand{\varepsint}[0]{(-\varepsilon, \varepsilon)}

\subjclass[2000]{
Primary 53C20, 
Secondary 
53C24, 
53C22, 
53B30, 
53C05, 
53C60, 
57R50} 

\newcommand{\npag}[0]{}

\newcommand{\slaz}[0]{\setminus \{0\}}

\title{Descending maps between slashed tangent bundles}
\keywords{geodesic conjugacy, geodesic flow, sprays, Finsler geometry, 
boundary rigidity, descending maps, complete lift, Jacobi fields}

\author[Bucataru]{Ioan Bucataru}
\address{Ioan Bucataru, Faculty of Mathematics, Al.I.Cuza University
B-dul Carol 11, Iasi, 700506, Romania}
\urladdr{http://www.math.uaic.ro/\textasciitilde{}bucataru/}

\author[Dahl]{Matias F. Dahl}
\address{Matias F. Dahl, Institute of Mathematics, P.O.Box 1100, 02015
Helsinki University of Technology, Finland}
\urladdr{http://www.math.tkk.fi/\textasciitilde{}fdahl/}

\date{\today}

\begin{document}
\begin{abstract}
  Suppose $TM\setminus \{0\}$ and $T\widetilde M\setminus\{0\}$ are
  slashed tangent bundles of two smooth manifolds $M$ and $\widetilde
  M$, respectively.  In this paper we characterize those
  diffeomorphisms $F\colon TM\setminus\{0\} \to T\widetilde
  M\setminus\{0\}$ that can be written as $F = (D\phi)|_{TM\slaz}$ for
  a diffeomorphism $\phi\colon M\to \wt M$. When $F =
  (D\phi)|_{TM\slaz}$ one say that $F$ \emph{descends}.  If $M$ is
  equipped with two sprays, we use the characterization to derive
  sufficient conditions that imply that $F$ descends to a totally
  geodesic map.  Specializing to Riemann geometry we also obtain
  sufficient conditions for $F$ to descent to an isometry.
\end{abstract}
\maketitle


\npag
\section{Introduction}
In this paper we study the following differential-topological problem: 
\begin{itemize}
\item[$(\ast)$]
  Suppose $M$ and $\widetilde M$ are smooth manifolds, and suppose
  that $\bigMap$ is a diffeomorphism between slashed tangent bundles
\begin{eqnarray}
\label{mainEq}
  \bigMap \colon TM\slaz &\to& T\widetilde M\slaz.
\end{eqnarray}
Characterize those maps $F$ that can be written as
$\bigMap=(D\phi)|_{TM\slaz}$ for a diffeomorphism $\phi\colon M\to
\widetilde M$, where $D\phi$ is the tangent map of $\phi$.
\end{itemize}
Problem $(\ast)$ is related to anisotropic boundary rigidity problems
on Riemannian manifolds \cite{Croke:2004, Uhlmann:2001, PeUh:2005}. It is also
the setting for studying conjugate geodesic flows. For an overview of
this topic for Riemann metrics, see \cite[p.~495]{Berger:panorama}.
When $\bigMap = (D\phi)|_{TM\slaz}$ one say that map $\bigMap$
\emph{descends} to a map $\phi\colon M\to \widetilde M$
\cite{Croke:2004}.  

Let us first note that if $f$ is a diffeomorphism between
(unslashed) cotangent bundles 
\begin{eqnarray*}
  f\colon T^\ast M &\to& T^\ast \widetilde M, 
\end{eqnarray*}
the analogous problem is well understood. Namely, $f$ descends into a
diffeomorphism $\psi\colon \widetilde M\to M$ if and only if $f$
preserves the canonical $1$-forms on $T^\ast M$ and $T^\ast \widetilde
M$, respectively. This result characterizes diffeomorphic
symplectomorphisms between cotangent bundles that arise from
diffeomorphisms between the base manifolds. The result can be seen as
a consequence of Euler's theorem for homogeneous
functions. Alternatively, $f$ defines a map between $M$ and $\wt M$
since $f$ maps zero covectors to zero covectors
\cite[p.~186]{AbrahamMarsden:1994}, \cite{Haro2000},
\cite[p.~66]{LibermannMarle:1987}, \cite[p.~212]{MiHrShSa}, and
\cite[p.~22]{daSilva:2000}. When $f$ is only defined between slashed
cotangent bundles this characterization is no longer valid
\cite[p.~34]{WeinsteinLec76}.

In this work we study maps $\bigMap$ as in equation
\eqref{mainEq}. Hence $\bigMap$ is defined and smooth only for
non-zero vectors. In this case the problem is more difficult
since we can not use the zero section to define a map $\phi\colon M\to
\widetilde M$. We can neither use Euler's theorem for homogeneous
functions to deduce that $\bigMap$ is linear in the vector variable.
Our first main result is Theorem \ref{lemmaKey}. It states that if
$\bigMap$ is a diffeomorphism $\bigMap\colon TM\slaz\to T\widetilde
M\slaz$, then $\bigMap=(D\phi)|_{TM\slaz}$ for a diffeomorphism
$\phi\colon M\to \wt M$ if and only if
\begin{eqnarray}
\label{DPhiAlgEq}
  D\bigMap &=& \wt \kappa_2\circ D\bigMap\circ \kappa_2,
\end{eqnarray}
where $\kappa_2$ and $\wt \kappa_2$ are the canonical involutions on
$TTM$ and $TT\wt M$ (Section \ref{sec:k2}). 

Let us note that Problem $(\ast)$ is a problem in
differential-topology.  Let us also note that Theorem \ref{lemmaKey}
provides a differential-topological answer.  One can interpret Theorem
\ref{lemmaKey} as an analogue to Poincar\'e's lemma for
diffeomorphisms; if the derivative of diffeomorphism $\bigMap$
satisfies algebraic condition \eqref{DPhiAlgEq}, then diffeomorphism
$\bigMap$ can be written as the derivative of another diffeomorphism.

As an application of Theorem \ref{lemmaKey} we prove 
Theorem
\ref{mainThmXXSpray} and Theorem \ref{mainThmXX}. In these theorems
we restate the assumption for a map to descend using mapping
properties for geometric objects of two sprays $S$ and $\wt S$ on $M$.
In Theorem \ref{mainThmXXSpray} we give sufficient
conditions that imply that $F$
descends to a totally geodesic map $\phi\colon M\to M$.  In Theorem
\ref{mainThmXX} we specialize to Riemann geometry and 
give sufficient conditions that imply that $F$
descends to an isometry $\phi\colon M\to M$.
The key assumption in both theorems is that $F$ maps Jacobi fields of
$S$ into Jacobi fields of $\wt S$.  This means that both theorems
essentially describe to what extend Jacobi fields determine the spray
(or Riemann metric).  Let us point out that Jacobi fields and
curvature are related. However, they are also different, since the
covariant derivative is needed to relate one to the other.  For
results on the relation between curvature and the Riemann metric, see
\cite{Kulkarni:1970}, \cite{Liu:1974}, \cite{Yau:1974}, and the
Cartan-Ambrose-Hicks theorem \cite[p.~31--34]{Cheeger:2008}. For the
real-analytic case, see also \cite[p.~259--261]{KobNom:1963} and
\cite{NomYan:1967}.

In Riemann geometry, the \emph{geodesic conjugacy problem} asks the
following \cite[p.~495]{Berger:panorama}: If $F\colon TM\slaz\to T
M\slaz$ maps integral curves of one Riemann metric into integral
curves of another Riemann metric, what additional assumptions are
required for $F$ to be induced by an isometry?  If $F$ satisfies the
assumptions in Theorem \ref{mainThmXX}, then $F$ necessarily maps
integral curves into integral curves (see Step $1$ in the proof of
Theorem \ref{mainThmXXSpray}).
Hence Theorem \ref{mainThmXX} is also a contribution to understanding
the geodesic conjugacy problem.

\section{Preliminaries}
\label{defsz}
By a manifold $M$ we mean a topological Hausdorff space with countable
base that is locally homeomorphic to $\setR^n$ with $C^\infty$-smooth
transition maps and $n=\dim M \ge 1$. All objects are assumed to be 
$C^\infty$-smooth where defined.

The next sections collect results about iterated tangent bundles
we will need. For a more detailed discussion and references we refer
to \cite{BucataruDahl:2008, BucataruDahl:2008:conjugate}.
 
\subsection{Iterated tangent bundles}
If $M$ is a manifold, let $TM$ be the tangent bundle of $M$.  For
$r\ge 0$, the $r$th \emph{iterated tangent bundle} $T^rM$ is defined
inductively by setting $T^rM=M$ when $r=0$, and $T^{r}M=T(T^{r-1}M)$
when $r\ge 1$.  Let $\pi_r$ be the canonical projection operators
$\pi_{r}\colon T^{r+1}M \to T^rM$ when $r \ge 0$.  Occasionally we
also write $\pi_{TTM\to M}$, $\pi_{TM\to M}, \ldots$ instead of
$\pi_0\circ \pi_1$, $\pi_0, \ldots$.  Unless otherwise specified, we
always use canonical local coordinates (induced by local coordinates
on $M$) for iterated tangent bundles.  If $x^i$ are local coordinates
for $M$, we denote induced local coordinates for $TM$, $TTM$, and
$TTTM$ by
\begin{eqnarray*}
& & (x,y), \\
& & (x,y,X,Y), \\
& & (x,y,X,Y,u,v,U,V).
\end{eqnarray*}
As above, we usually leave out indices for local coordinates and write
$(x,y)$ instead of $(x^i, y^i)$. For $p\in M$ let $T_pM =
\pi_0^{-1}(p)$.

For $r\ge 1$, we treat $T^rM$ as a vector bundle over the manifold
$T^{r-1}M$ with the vector space structure induced by projection
$\pi_{r-1}\colon T^rM\to T^{r-1}M$. Thus, if $\{ x^i : i=1, \ldots,
2^{r-1}n \}$ are local coordinates for $T^{r-1}M$, and $(x,y)$ are
local coordinates for $T^rM$, then vector addition and scalar
multiplication are given by
\begin{eqnarray}
\label{eq:TMplus}
  (x,y) + (x, \widetilde y) &=& (x,y+\widetilde y), \\
\label{eq:mult}
  \lambda \cdot (x,y)  &=& (x,\lambda  y).
\end{eqnarray}
For $r\ge 0$, a \emph{vector field} on an open set $B\subset T^rM$
is a smooth map $X\colon B\to T^{r+1}M$ such that $\pi_{r}\circ X =
\operatorname{id}_{B}$.  The set of all vector fields on $B$ is
denoted by $\vfield{B}$.
Suppose that $\gamma$ is a smooth map $\gamma\colon
\left(-\varepsilon, \varepsilon\right)^k \rightarrow T^rM $ where
$k\ge 1$ and $r\ge 0$. If $\gamma(t^1,\ldots, t^k) =
(z^i(t^1, \ldots, t^k ))$ in local coordinates $(z^i)$ for $T^rM$,
then the \emph{derivative} of $\gamma$ with respect to
variable $t^j$ is the map $\partial_{t^j}\gamma \colon
\left(-\varepsilon, \varepsilon\right)^k$ $\to T^{r+1}M$ defined by
$\partial_{t^j}\gamma=\left(z^i, {\partial z^i}/{\partial
t^j}\right)$. When $k=1$ we also write $\gamma'=\partial_{t}\gamma$
and say that $\gamma'$ is the \emph{tangent of $\gamma$}.
If $f\colon T^rM\to T^s\wt M$ ($r,s\ge 0$) is a map between 
iterated tangent bundles
and $c\colon I\to T^r M$ is a curve, then 
\begin{eqnarray}
\label{eq:chainRule}
  (f\circ c)'(t) &=& Df \circ c'(t), \quad t\in I.
\end{eqnarray}
Unless otherwise stated we always assume that $I$ is an open interval
in $\setR$ (and we do not exclude unbounded intervals).

If $\xi\in T^rM$ for $r\ge 2$, then there exists
a  map $V\colon \varepsint^2\to T^{r-2}M$ such that
\begin{eqnarray}
\label{xi:expression}
  \xi = \partial_t \partial_s V(t,s) |_{t=s=0}.
\end{eqnarray}

\subsection{Canonical involution}
\label{sec:k2}
On the iterated tangent bundle $T^rM$ where $r\ge 2$ the \emph{canonical
  involution} is the unique diffeomorphism $\kappa_r \colon T^rM\to
T^rM$ such that
\begin{eqnarray}
\label{kr_comm_rule}
  \partial_s   \partial_t c(t,s) &=& \kappa_r  \circ \partial_t   \partial_s c(t,s)
\end{eqnarray}
for all smooth maps $c\colon \varepsint^2\to T^{r-2}M$. Let also
$\kappa_1=\operatorname{id}_{TM}$.  
In local coordinates for $TTM$ and $TTTM$, it follows that
\begin{eqnarray*}
  \kappa_2(x,y,X,Y) &=& (x,X,y,Y), \\
  \kappa_3(x,y,X,Y,u,v,U,V) &=& (x,y,u,v,X,Y,U,V).
\end{eqnarray*}
For any $r\ge 1$, we have
\begin{eqnarray}
\label{commRel}
  \pi_{r} \circ \kappa_{r+1} &=& D\pi_{r-1}, \\
\label{commRelA}
  \pi_{r-1} \circ \pi_r\circ \kappa_{r+1} &=& \pi_{r-1}\circ \pi_r.
\end{eqnarray}
If $\phi$ is a map $\phi\colon M\to \wt M$, then equations
\eqref{eq:chainRule}, \eqref{xi:expression}, and \eqref{kr_comm_rule} imply that
\begin{eqnarray}
\label{kapComA}
\wt \kappa_{2}\circ DD\phi\circ \kappa_{2} &=& DD\phi.
\end{eqnarray}
As in equation \eqref{kapComA} we denote
involution operators on $T^rM$ and $T^r\wt M$ by $\kappa_r$ and $\wt
\kappa_r$, respectively. Similarly, we denote projection operators
by $\pi_r$ and $\wt \pi_r$. 

\subsection{Slashed tangent bundles}
The \emph{slashed tangent bundle} for $M$ is defined as the open set of
non-zero vectors, 
\begin{eqnarray*}
  TM\slaz &=& \{ \xi \in TM : \xi\neq 0\}.
\end{eqnarray*}
For $r\ge 2$ we generalize and define
\begin{eqnarray*}
  T^rM\slaz &=& \{ \xi \in T^rM : (D\pi_{T^{r-1}M\to M})(\xi) \in TM\slaz \}.
\end{eqnarray*}
When $r\ge 2$, $\kappa_r$ restricts to a diffeomorphism
\begin{eqnarray}
\label{eq:kadi}
  \kappa_r\colon T^rM\slaz &\to& T(T^{r-1}M\slaz).
\end{eqnarray}
If $\bigMap$ is a map $\bigMap\colon TM\slaz \to T\wt M\slaz$,
then
\begin{eqnarray}
\label{kapComBB}
\wt \kappa_{3}\circ DD\bigMap\circ \kappa_{3} &=& DD\bigMap
\ \  \mbox{on} \ \ TT(TM\slaz).
\end{eqnarray}

\section{A differential-topological characterization}
\label{charDiff}
Theorem \ref{lemmaKey} is the first main result in this paper. The
theorem is a differential-topological characterization of descending
maps between slashed tangent bundles.
\begin{theorem}
\label{lemmaKey}
Suppose $M$ and $\widetilde M$ are smooth manifolds.
If $\bigMap$ is a smooth map $\bigMap\colon TM\slaz\to T\widetilde
M\slaz$, then the following conditions are equivalent:
\begin{enumerate}
\item There exists a smooth map $\phi\colon M\to \widetilde M$ such
  that
\begin{eqnarray*}
\label{eq:desDphi}
  \bigMap &=& (D\phi)|_{TM\slaz}.
\end{eqnarray*}
\item On $TTM\slaz\cap T(TM\slaz)$,
\begin{eqnarray*}
\label{eq:poincare} 
  D\bigMap &=& \wt \kappa_2\circ D\bigMap\circ \kappa_2.
\end{eqnarray*}
\end{enumerate}
What is more, if $\bigMap$ is a diffeomorphism, and $\phi$ exists, then $\phi$
is a diffeomorphism.
\end{theorem}  

Let us make three remarks about Theorem \ref{lemmaKey} assuming that
$\phi$ exists. First, when $\phi$ exists, it is unique, and the
following diagram commutes:
\begin{eqnarray*}
\begin{xy}
\xymatrix{
  TM\slaz\ar@{->}[rr]^{\bigMap}\ar[d]_{\pi_0} & & T\widetilde M\slaz \ar[d]^{\wt \pi_0} \\
  M\ar@{->}[rr]_{\phi} & &\wt M
}
\end{xy}
\end{eqnarray*}
Second, since $F$ is a map between slashed tangent bundles, $\phi$ is
necessarily an immersion. Thus, if $\dim M = \dim \wt M$, the inverse
function theorem implies that $\phi$ is a local diffeomorphism.
Third, if $\phi$ is a diffeomorphism, then equation $\bigMap = D\phi$
extends $\bigMap$ into a (smooth) diffeomorphism $\bigMap\colon TM\to
T\wt M$.

Theorem \ref{lemmaKey} is a direct consequence of the next two lemmas;
implication \emii $\Rightarrow$ \emi follows by Lemma \ref{lemma:desExists},
the last claim follows by Lemma \ref{secDescent}, and 
the easy implication \emi $\Rightarrow$ \emii follows by equation
\eqref{kapComA}.

\begin{lemma}
\label{lemma:desExists}
Let $\bigMap$ be a smooth map $\bigMap\colon
TM\slaz\to T\widetilde M\slaz$ that satisfies condition (ii) in
Theorem \ref{lemmaKey}, and let $\phi$ be the set-valued map $\phi\colon
M\to P(\wt M)$,
\begin{eqnarray*}
  \phi(p) &=& \widetilde \pi_0 \circ \bigMap(T_pM\slaz), \quad p\in M,
\end{eqnarray*}
where $P(\wt M)$ is the power set of $\wt M$. Then
\begin{enumerate}
\item $\phi$ defines a smooth single-valued map $\phi\colon M\to \wt M$,
\item $\bigMap= (D\phi)|_{TM\slaz}$.
\end{enumerate}
\end{lemma}

\begin{proof}
  To show that $\phi$ is single-valued we show that map $C\colon
  T_pM\slaz \to \wt M$,
\begin{eqnarray*}
  C(\xi) &=& \widetilde \pi_0 \circ \bigMap(\xi), \quad \xi\in T_pM\slaz,
\end{eqnarray*}
is constant when $p\in M$ is fixed. If $\xi, \eta\in T_pM\slaz$ we can
find a $w\in TTM\slaz\cap T(TM\slaz)$ such that $\pi_{1}(w) = \xi$ and
$D\pi_{0}(w) = \eta$.  Using equations 
\eqref{commRel} and \eqref{commRelA}, and the
assumption on $D\bigMap$ we have
\begin{eqnarray*}
  C(\xi) &=& \wt \pi_0\circ \bigMap\circ \pi_1(w) \\
 &=& \wt \pi_0 \circ \wt \pi_1 \circ D\bigMap (w) \\
 &=& \wt \pi_0 \circ \wt \pi_1 \circ \wt \kappa_2\circ D\bigMap\circ\kappa_2 (w) \\
 &=& \wt \pi_0 \circ \wt \pi_1 \circ D\bigMap\circ\kappa_2 (w) \\
 &=& \wt \pi_0 \circ \bigMap \circ \pi_1 \circ\kappa_2 (w) \\
 &=& C(\eta),
\end{eqnarray*}
and $\phi$ defines a single-valued map $\phi\colon M\to \wt M$.  If
$p\in M$, and $U$ is a non-vanishing vector field $U\in
\vfield{B}$ defined in a neighborhood $B\subset M$ of $p$, then
\begin{eqnarray*}
  \phi(x) &=& \widetilde \pi_0 \circ \bigMap\circ U(x),\quad x\in B,
\end{eqnarray*}
and $\phi$ is smooth near $p$. 
To prove \emph{(ii)}, let $\xi\in T_pM\slaz$, and let $U$ be a
non-vanishing vector field defined near $p$ such that $U(p)=\xi$.
Starting from $D\phi(\xi) = D(\wt \pi_0 \circ F\circ U)(\xi)$, a
similar calculation\proofread{ Using equations
  \eqref{commRel}, \eqref{commRelA}, and the assumption on $D\bigMap$
  we have
\begin{eqnarray*}
  D\phi(\xi) &=& D\wt \pi_0 \circ D\bigMap(w) \\
             &=& D\wt \pi_0 \circ \wt \kappa_2\circ D\bigMap\circ \kappa_2 (w) \\
             &=& \wt \pi_1 \circ D\bigMap\circ \kappa_2 (w) \\
             &=& \bigMap\circ \pi_1\circ \kappa_2 (w) \\
             &=& \bigMap\circ D\pi_0 (w) \\
             &=& \bigMap(\xi).
\end{eqnarray*}
} used to prove that map $C$ is constant shows that
$D\phi(\xi)=\bigMap(\xi)$.
\end{proof}

\begin{lemma} 
\label{secDescent}
If $\bigMap\colon TM\slaz\to T\widetilde M\slaz$ is a
diffeomorphism, and $\bigMap= (D\phi)|_{TM\slaz}$ for a smooth map
$\phi\colon M\to \widetilde M$, then $\phi$ is a diffeomorphism.
\end{lemma}

\begin{proof}
  Since $F$ is a diffeomorphism, we have $\dim M = \dim \wt M$, and by
  the inverse function theorem, $\phi$ is a local diffeomorphism.
  If $\xi \in TT\wt M\slaz\cap T(T\wt M\slaz)$, then there exists
  a $\zeta \in T(TM\slaz)$, such that $\xi = DF(\zeta)$. 
  If $\zeta = \gamma'(0)$ for a curve $\gamma\colon \varepsint\to
  TM\slaz$, we obtain $0 \neq D\pi_0(\xi)= D\phi\circ
  D\pi_0(\zeta)$. Hence $\zeta \in TTM\slaz \cap T(TM\slaz)$, so
  $DF(\zeta)=\wt \kappa_2 \circ DF\circ \kappa_2(\zeta)$, and
\begin{eqnarray*}
  \kappa_2\circ D(F^{-1})\circ \wt \kappa_2(\xi) &=& D(F^{-1})(\xi).
\end{eqnarray*}
By Lemma \ref{lemma:desExists}, there exists a smooth map $\rho\colon
\wt M\to M$ such that $F^{-1} = D\rho|_{T\wt M\slaz}$. Since $\rho
\circ \phi=\operatorname{id}|_M$ and $\phi \circ
\rho=\operatorname{id}|_{\wt M}$, it follows that $\phi$ is a
diffeomorphism.
\end{proof}

\section{Sprays}
\label{sec:sprays}
The motivation for studying sprays is that they provide a unified
framework for studying geodesics for Riemannian metrics, Finsler
metrics, and non-linear connections. See \cite{BucataruMiron:2007,
  Sakai1992, Shen2001}. Following \cite{BucataruDahl:2008,
  BucataruDahl:2008:conjugate} we next define a spray on an iterated
tangent bundle $T^rM$.

\begin{definition}[Spray]
\label{def:spray}
A \emph{spray} on $T^rM$ where $r\ge 0$ is a vector field $S\in
\vfield{T^{r+1}M\slaz}$ such that $\kappa_{r+2}\circ S = S$ and
$[S,\setC_{r+1}]=S$, where $\setC_r\in \vfield{T^rM}$, $r\ge 1$ is the
\emph{Liouville vector field} defined by
\begin{eqnarray*}
  \setC_r(\xi) &=& \partial_t (\xi + t \xi)|_{t=0}, \quad \xi\in T^rM.
\end{eqnarray*}
\end{definition}

If $(x,y,X,Y)$ are local coordinates for $T^{r+2}M$ then a spray $S$
can be written as 
\begin{eqnarray}
\nonumber
  S(x,y) &=& (x,y,y,-2G^i(x,y)) \\
\label{eq:geo}
         &=& y^i\pd{}{x^i} - 2G^i(x,y) \pd{}{y^i}
\end{eqnarray}
for locally defined component functions $G^i\colon T^{r+1}M\slaz\to \setR$
that are positively $2$-homogeneous. That is, 
\begin{eqnarray*}
  G^i(\lambda  y) &=& \lambda^2  G^i(y), \quad \lambda>0.
\end{eqnarray*}
A curve $c\colon I\to T^rM$ is \emph{regular} if $c'(t)
\in T^{r+1}M\slaz$ for all $t\in I$. 
That is, curve $c$ is regular if and only
if its projection $\pi_{T^rM\to M}\circ c\colon I\to M$ is regular.

\begin{definition}[Geodesic]
\label{def:geo}
If $S$ is a spray on $T^rM$ for $r\ge 0$, a regular curve $c\colon
I\to T^rM$ is a \emph{geodesic} if 
\begin{eqnarray*}
  c'' &=& S\circ c'.
\end{eqnarray*}
\end{definition}

That is, a regular curve $c$ is a geodesic of spray $S$ if and only if
$c'$ is an integral curve of $S$. Conversely, suppose that
$\gamma\colon I\to T^{r+1}M\slaz$ is an integral curve of $S$, whence
$\gamma'=S\circ \gamma$. Since $\kappa_{r+2}\circ S=S$, there
is\proofread{ $c'=(\pi_r\circ \gamma)'=D\pi_r\circ \gamma' =
  \pi_{r+1}\circ \kappa_{r+2}\circ \gamma' = \pi_{r+1}\circ
  \kappa_{r+2}\circ S\circ \gamma=\pi_{r+1}\circ S\circ \gamma =
  \gamma$.}
a geodesic $c\colon I\to T^{r}M$, $c=\pi_{r}\circ \gamma$ such that
$\gamma=c'$.

Any geodesic $c\colon I\to T^rM$ of a spray $S$ is uniquely
determined by $c'(t_0)$ for one $t_0\in I$. The \emph{geodesic flow}
of a spray $S$ is defined as the flow of $S$ as a vector field, and a
spray is \emph{complete} if $S$ is complete as a vector field.

If $S$ is locally written as in equation \eqref{eq:geo} and $c(t)=(x^i(t))$, 
then $c$ is a geodesic if and only if 
\begin{eqnarray*}
  \ddot x^i(t)  + 2G^i\circ c'(t) &=& 0.
\end{eqnarray*}

\subsection{Jacobi fields}
We define Jacobi fields for a spray using the complete lift 
following \cite{BucataruDahl:2008, BucataruDahl:2008:conjugate}.
See also \cite{Lewis:2001:GeomMax, Michor1996, Yano1973}. 

\begin{definition}[Complete lift] 
  The \emph{complete lift} of a spray $S$ on $M$ is the spray $S^c\in
  \vfield{TTM\slaz}$ on $TM$ given by
\begin{eqnarray}
\label{eq:Sc}
  S^c &=& D\kappa_{2}\circ \kappa_{3}\circ DS\circ \kappa_{2}.
\end{eqnarray}
\end{definition}
Suppose that $S$ is locally given by equation \eqref{eq:geo}. Then
$S^c$ is locally given by
\begin{eqnarray*}
\label{compSc}
  S^c &=& \left(x,y,X,Y,X,Y,-2A^i(x,y,X,Y), -2 B^i(x,y,X,Y) \right) \\
\nonumber
      &=& X^i \pd{}{x^i} + Y^i \pd{}{y^i}  -2 A^i(x,y,X,Y) \pd{}{X^i} -2 B^i(x,y,X,Y)  \pd{}{Y^i},
\end{eqnarray*}
where $A^i$ and $B^i$ are vertical and complete lifts of functions $G^i$
\cite{BucataruDahl:2008:conjugate}, 
\begin{eqnarray*}
  A^i(x,y,X,Y) &=& G^i(x,X), \\  
  B^i(x,y,X,Y) &=& \pd{G^i}{x^a}(x,X)y^a + \pd{G^i}{y^a}(x,X)Y^a.
\end{eqnarray*}
Spray $S^c$ is complete if and only if spray $S$ is complete.

\begin{definition}[Jacobi field] 
\label{def:Jacobi}
Suppose $S$ is a spray on $M$.
A \emph{Jacobi field} for $S$ is a geodesic $J\colon I\to TM$ of $S^c$. 
\end{definition}

If $J\colon I\to TM$ is a Jacobi field for $S$, then 
curve $c\colon I\to M$, $c=\pi_0\circ J$ is a geodesic for $S$ and
we say that \emph{$J$ is a Jacobi field along $c$}.
Next we show that Definition \ref{def:Jacobi}
coincides with the usual characterization of Jacobi fields in terms of
geodesic variations. For proofs and 
discussions, see
\cite{BucataruDahl:2008, BucataruDahl:2008:conjugate}.

\begin{definition}[Geodesic variation] 
  Suppose $S$ is a spray on $M$, and $c \colon I \to M$ is a
 geodesic for $S$.  Then a \emph{geodesic variation} of $c$ is a
 smooth map $V\colon I\times \varepsint\to M$ such that
\begin{enumerate}
\item $V(t,0)=c(t)$ for all $t\in I$, 
\item $t\mapsto V(t,s)$ is a geodesic for all $s\in \varepsint$.
\end{enumerate}
\end{definition}

Suppose that $I$ is a closed interval. 
Then we say that a curve $J\colon I\to TM$ is a Jacobi field if we
can extend $J$ into a Jacobi field defined on an open interval.
Similarly, a map $V\colon I\times \varepsint\to M$ is a geodesic
variation if there is a geodesic variation $V^\ast\colon I^\ast \times
(-\varepsilon^\ast, \varepsilon^\ast)\to M$ such that $V=V^\ast$ on
the common domain of $V$ and $V^\ast$ and $I\subset I^\ast$.

\begin{proposition}[Jacobi fields and geodesic variations]
\label{CharacterizationOfJacobiFields}
Let $S$ be a spray on $M$, let $J\colon I\to
TM$ be a curve, where $I$ is open or closed, and let
$c\colon I\to M$ be the curve $c = \pi_0 \circ J$.
\begin{enumerate}
\item If $J$ can be written as
\begin{eqnarray}
\label{eq:JacDef}
  J(t) &=& \left.\partial_s V(t,s)\right|_{s=0}, \quad t\in I
\end{eqnarray}
for a geodesic variation $V\colon I\times \varepsint\to M$, 
then $J$ is a Jacobi field along $c$.
\label{CharacterizationOfJacobiFields:Jac2Sc}
\item If $I$ is compact and $J$ is a Jacobi field along $c$, then
      there exists a geodesic variation $V\colon I\times \varepsint\to
      M$ such that equation \eqref{eq:JacDef} holds.
\label{CharacterizationOfJacobiFields:Sc2Jac}
\end{enumerate}
\end{proposition}

\begin{remark}[Zero Jacobi field]
\label{zeroJacFieldRem}
If $c\colon I\to M$ is a geodesic for a spray $S$, then the \emph{zero
  Jacobi field} along $c$ is the Jacobi field $J\colon I\to TM$ that
is locally induced by the constant geodesic variation $V(t,s)=c(t)$. 
Globally,
\begin{eqnarray*}
  J(t) &=& D\pi_0 \circ \setC_1 \circ c'(t), \quad t\in I.
\end{eqnarray*}
If zeroes of a Jacobi fields converge, then the Jacobi field is a zero
Jacobi field.
\end{remark}

\npag
\section{Maps that preserve structure}
\label{maps}
Throughout this section we assume that $S$ and $\widetilde S$ are
sprays on manifolds $M$ and $\widetilde M$, respectively. We proceed
by studying maps 
that preserve \emph{(i)} integral curves, \emph{(ii)} geodesics, and
\emph{(iii)} Jacobi fields. In Section \ref{sec:Riemann} we will also
study maps between Riemann manifolds that preserve inner products.

\subsection{Maps that preserve integral curves}
\label{sec:conjugateMaps}
We say that a map
\begin{eqnarray*}
  \bigMap\colon TM\slaz &\to& T\widetilde M\slaz
\end{eqnarray*}
\emph{preserves integral curves} if $\bigMap\circ \gamma\colon I\to
T\widetilde M\slaz$ is an integral curve of $\widetilde S$ whenever
$\gamma \colon I\to TM\slaz$ is an integral curve of $S$. When such a
map $\bigMap$ exists, we say that sprays $S$ and $\widetilde S$ are
\emph{conjugate}. Condition \emiii in the next proposition shows that
this corresponds to the usual definition of geodesic conjugacy in Riemann
geometry \cite{Berger:panorama, Croke:2004, Uhlmann:2001}.

\begin{proposition}
\label{prop:IntPres}
Suppose $\bigMap$ is a smooth map $\bigMap\colon TM\slaz \to T\widetilde M\slaz$.
Then the following conditions are equivalent:\proofread{
Locally the condition reads:
\begin{enumerate}
\item[(iv)] If $S$ and $\widetilde S$ have local expressions
\begin{eqnarray*}
  S &=& (x^i,y^i,y^i,-2G^i(x,y)), \quad \widetilde S = (\widetilde x^i,\widetilde y^i,\widetilde y^i,-2\widetilde G^i(\widetilde x,\widetilde y))
\end{eqnarray*}
and locally $\bigMap(x,y) = (F_1^i(x,y), F_2^i(x,y))$, then
\begin{eqnarray*}
   F_2^i(x,y) &=& \pd{F_1^i}{x^r}(x,y) y^r -2 \pd{F_1^i}{y^r}(x,y) G^r(x,y), \\
   \widetilde G^i\circ \bigMap(x,y) &=& \pd{F_2^i}{y^r}(x,y) G^r(x,y) - \frac{1}{2} \pd{F_2^i}{x^r}(x,y) y^r.
\end{eqnarray*}
\end{enumerate}
}
\begin{enumerate}
\item $\bigMap$ preserves integral curves.
\item $\widetilde S\circ \bigMap = D\bigMap \circ S$ on $TM\slaz$.
\item If $\Phi_t$ and $\widetilde \Phi_t$ are geodesic flows of $S$ and
$\widetilde S$, respectively, then 
the following diagram commutes: 
\begin{eqnarray*}
\begin{xy}
\xymatrix{
  TM\slaz\ar@{->}[r]^{\bigMap}\ar[d]_{\Phi_t} & T\widetilde M\slaz \ar[d]^{\widetilde \Phi_t} \\
  TM\slaz\ar@{->}[r]_{\bigMap} & T\widetilde M\slaz
}
\end{xy}
\end{eqnarray*}
\end{enumerate}
\end{proposition}

\subsection{Maps that preserve geodesics} 
We say that a map 
\begin{eqnarray*}
  \phi\colon M &\to& \widetilde M
\end{eqnarray*} 
is a \emph{totally geodesic map} if $\phi\circ c\colon I\to \widetilde
M$ is a geodesic for $\widetilde S$ whenever $c\colon I\to M$ is a
geodesic for $S$ \cite[Chapter 6]{KobNom:1963}.

In Definition \ref{def:geo}, we assume that geodesics are regular
curves.  If $\phi$ is a totally geodesic map, we can therefore
restrict $D\phi$ to a map $D\phi\colon TM\slaz\to T\wt M\slaz$. Hence
every totally geodesic map $\phi$ is an immersion, and if $\dim M =
\dim \wt M$, then $\phi$ is also a local diffeomorphism.
The definition of a totally geodesic map does not depend on
derivatives of $\phi$. However, if $\phi\colon M\to \widetilde M$ is a
homeomorphism, it follows that $\phi$ is a diffeomorphism \cite{Brickell}.  

\begin{proposition} 
\label{prop:TGIN}
Suppose $\phi\colon M\to \widetilde M$ is a smooth immersion.  Then
$\phi$ is a totally geodesic map if and only if restriction
$D\phi\colon TM\slaz\to T\wt M\slaz$ preserves integral curves.
\end{proposition}

\subsection{Maps that preserve Jacobi fields}
\label{sec:StrJacobi}
We say that a map 
\begin{eqnarray*}
  \bigMap\colon TM\slaz &\to& T\widetilde M\slaz
\end{eqnarray*}
\emph{preserves Jacobi fields} if for any Jacobi field
$J\colon I \to TM\slaz$ without zeroes, 
\begin{eqnarray}
\label{phiDefEq}
  \widetilde J(t) &=& \bigMap\circ J(t), \quad t\in I
\end{eqnarray}
is a Jacobi field $\widetilde J\colon I \to T\widetilde
M\slaz$ without zeroes.

In the above definition, we only apply $\bigMap$ to Jacobi fields
without zeroes.
The next proposition shows that we can still map Jacobi fields with
isolated zeroes.

\begin{proposition} 
\label{preJacobi_I}
Suppose $\wt S$ is complete, $\dim M\ge 2$, and $\bigMap$ is a map
$\bigMap\colon TM\slaz \to T\widetilde M\slaz$ that preserves Jacobi
fields. If $J\colon \setR \to TM$ is a Jacobi field for $S$ that is not
identically zero, then there exists a Jacobi field $\widetilde J\colon
\setR\to T\widetilde M$ for $\wt S$ such that
\begin{eqnarray}
\label{JtildeExists}
  \widetilde J'(t) &=& D\bigMap\circ J'(t), \quad t\in \setR\setminus Z,
\end{eqnarray}
where $Z=\{ t\in \setR: J(t)=0\}$.
\end{proposition}

The proof of Proposition \ref{preJacobi_I} is slightly technical and
is given in Appendix \ref{appZero}. The idea of the proof is to
approximate a Jacobi field $J$ with an isolated zero by a variation of
Jacobi fields without zeroes (see Lemma
\ref{lemma:puncturedJacobiVariation}). Then $\bigMap$ maps each non-zero
Jacobi field in the variation into a non-zero Jacobi field, and a
continuity argument shows that there exists a Jacobi field $\wt J$ as
in equation \eqref{JtildeExists}.

\begin{proposition} 
\label{preJacobi_II}
Suppose that map $\bigMap\colon TM\slaz \to T\widetilde M\slaz$ preserves
integral curves, and suppose that $J\colon I\to TM$ is a Jacobi
field for $S$. Then curve $\widetilde J\colon I\to T\widetilde M$,
\begin{eqnarray*}
   \widetilde J'(t) &=& \wt \kappa_2\circ D\bigMap\circ \kappa_2\circ J'(t), \quad t\in I
\end{eqnarray*} 
is a Jacobi field for $\wt S$.
\end{proposition}

\begin{proof} 
  Equation \eqref{eq:kadi} shows that curve $\wt J'\colon I\to TT\wt M
  \slaz$ is smooth. Proposition \ref{prop:IntPres} and equations
  \eqref{eq:kadi}, \eqref{kapComBB} and \eqref{eq:Sc} imply that
\begin{eqnarray*}
  \widetilde S^c \circ (\wt\kappa_2\circ D\bigMap \circ \kappa_2) &=& D(\wt \kappa_2\circ D\bigMap\circ \kappa_2) \circ S^c \ \mbox{on} \ TTM\slaz,
\end{eqnarray*}
and $\wt\kappa_2\circ D\bigMap \circ \kappa_2$
maps integral curves of $S^c$ into integral curves of $\wt S^c$.
\end{proof}

The next proposition is analogous to Proposition \ref{prop:TGIN}.

\begin{proposition} 
\label{DphiandtotGeo}
Suppose $\phi\colon M\to \widetilde M$ is a smooth immersion. Then $\phi$ is
a totally geodesic map if and only if restriction $D\phi\colon TM \slaz \to
T\widetilde M\slaz$ preserves Jacobi fields. 
 \end{proposition}

\begin{proof}
  If $\phi$ is totally geodesic, then Propositions \ref{prop:TGIN} and
  \ref{preJacobi_II} imply that $D\phi$ preserves Jacobi fields.
  For the converse direction, suppose that $D\phi$ preserves Jacobi
  fields and $c\colon I\to M$ is a geodesic for $S$. Then $c'$ is a
  Jacobi field for $S$, so $(D\phi)\circ c'$ is a Jacobi field for
  $\widetilde S$, and $\wt c = \wt \pi_0 \circ ( \phi\circ
  c)'=\phi\circ c$ 
is a geodesic $\wt c\colon I\to \wt M$ for
  $\widetilde S$.
\end{proof}

\section{Descending maps for sprays}
\label{descDescTrap}
In this section we prove Theorem \ref{mainThmXXSpray}, which gives 
sufficient conditions for a map $F\colon TM\slaz\to TM\slaz$ to
descend to a totally geodesic map between two sprays.  To formulate
the assumptions in Theorem \ref{mainThmXXSpray} we need the concept of
a \emph{trapping hypersurface}. This term is adapted from the concept
of a non-trapping manifold with boundary.
\begin{definition}[Trapping hypersurface] Suppose $S$ is a spray on a
  manifold $M$. A hypersurface $\Sigma \subset M$ is a \emph{trapping
    hypersurface for $S$} if for any $y\in TM\slaz$ there
  exists a geodesic $c\colon I \to M$ such that $c'(0) = y$ and
  $c(t)\in \Sigma$ for some $t\in I$. 
\end{definition}
The existence of a trapping hypersurface $\Sigma$ imposes a global
restriction on the behavior of geodesics. Namely, every geodesic
has to intersect $\Sigma$. An interpretation is that if
geodesics describe propagation of light, then the whole manifold is
visible from the trapping hypersurface. 

One way to construct a spray with a trapping hypersurfaces one can
start with two sprays on a manifold $B$ with boundary $\partial
B$. Using a \emph{smooth double} one can glue together two copies of
$B$ by identifying their boundary points. This gives a smooth manifold
$M$ without boundary that contains two copies of the interior of $B$
and one copy of boundary $\partial B$. See \cite[p.~184]{Hirsch:1976},
\cite[p.~463]{Lee:abc}, or \cite[p.~39]{Matsumoto:2002}. Assuming that
the two sprays are \emph{non­trapping} (see \cite{Dairbekov2006} for
the Riemann case), and assuming that they satisfy suitable
compatibility conditions on the boundary, one can glue together the
sprays into a spray on $M$ such that boundary $\partial B\subset M$ is
a trapping hypersurface.  For example, any great circle on the
$2$-sphere with the induced Euclidean metric is a trapping
hypersurface.

\begin{theorem}
\label{mainThmXXSpray}
Suppose $S$ and $\widetilde S$ are complete sprays on a manifold $M$
with $\dim M$ $\ge 2$. Furthermore, suppose that there exists a smooth
map $\bigMap\colon TM\slaz\to TM\slaz$ and a trapping hypersurface
$\Sigma\subset M$ for $S$ such that
\begin{enumerate}
\item 
$\bigMap$ maps Jacobi fields for
  $S$ into Jacobi fields for $\widetilde S$ (see Section \ref{sec:StrJacobi}),
\label{pr:jacobi}
\item 
for all $p\in \Sigma$, 
\begin{eqnarray}
\label{pr:SeqStil_xx}
  S(y) &=& \wt S(y), \quad\ \  y\in T_pM\slaz, \\
\label{pr:DFid_xx}
  D\bigMap(\xi) &=& \xi, \quad \quad\quad\xi\in T(T_pM\slaz).
\end{eqnarray}
\end{enumerate}
Then there exists a smooth map $\phi\colon M\to M$ such that
\begin{enumerate}
\item $\bigMap= (D\phi)|_{TM\slaz}$,
\item $\phi$ is a totally geodesic map (that maps geodesics for $S$
  into geodesics for $\widetilde S$).
\end{enumerate}
What is more, if $\bigMap$ is a diffeomorphism, then $\phi$ is a
diffeomorphism.
\end{theorem}

In the proof below, Subcase B also proves Subcase A. However, Subcase
A is included as it illustrates the main argument with minimal
technical detail.

\begin{proof} The proof is divided into two steps:

\noindent \emph{Step 1: Map $F$ maps integral curves of $S$
  into integral curves of $\wt S$.}\\
\noindent Let $c'\colon \setR\to TM\slaz$ be an integral curve of $S$,
where $c$ is a geodesic $c\colon \setR\to M$ of $S$. Then $c'$ is a
non-zero Jacobi field for $S$, and by assumption \ref{pr:jacobi},
$J=F\circ c'$, $J\colon \setR\to TM\slaz$ is a Jacobi field of $\wt S$
without zeroes. Since $\Sigma$ is trapping, there exists a $t_0\in
\setR$ such that $c(t_0)\in \Sigma$. By equations \eqref{eq:chainRule}
and \eqref{pr:DFid_xx}, we have $J'(t_0) = c''(t_0)$. If $\eta\colon
\setR\to TM\slaz$ is the integral curve of $\wt S$ determined by
$\eta(t_0) = c'(t_0)$, then $ \eta'(t_0) = \wt S\circ \eta(t_0) =
J'(t_0) $ by equation \eqref{pr:SeqStil_xx}.  Thus Jacobi fields
$\eta$ and $J$ coincide and $J$ is an integral curve of $\wt S$.

\noindent\emph{Step 2: If $\xi\in TTM\slaz\cap T(TM\slaz)$ we claim that}
\begin{eqnarray}
\label{eq:toProve}
  D\bigMap(\xi) &=& \kappa_2\circ D\bigMap\circ \kappa_2(\xi).
\end{eqnarray}
If equation \eqref{eq:toProve} holds, Theorem \ref{lemmaKey} implies
that $\bigMap = (D\phi)|_{TM\slaz}$ for a map $\phi\colon M\to M$,
whence assumption \ref{pr:jacobi} and Proposition \ref{DphiandtotGeo}
imply that $\phi$ is totally geodesic. (Alternatively, one could use
Step $1$ and Proposition \ref{prop:TGIN}.)  The last claim follows by
Theorem \ref{lemmaKey}. 

To prove equation \eqref{eq:toProve}, let $J\colon \setR\to TM$ be the
Jacobi field with $J'(0)= \xi$, and let $c\colon \setR\to M$ be
the geodesic $c=\pi_0 \circ J$. Since $\Sigma$ is a trapping
hypersurface, there is a $t_0\in \setR$ such that $c(t_0)\in \Sigma$.

\noindent\emph{Subcase A:} $J(t_0)\neq 0$.\\
Propositions \ref{preJacobi_I} and
\ref{preJacobi_II} imply that there exist Jacobi fields $J_1, J_2 \colon \setR\to TM$ for $\widetilde S$ such that
\begin{eqnarray*}
  J'_1(t) &=& D\bigMap\circ J'(t), \quad \quad\quad\quad\quad\mbox{when}\, t\in \setR \ \mbox{and} \  J(t)\neq 0, \\
  J'_2(t) &=& \kappa_2\circ D\bigMap\circ \kappa_2 \circ J'(t), \quad \mbox{when}\, t\in \setR.
\end{eqnarray*}
Since $J(t_0)\neq 0$, we have
\begin{eqnarray*}
  J'(t_0) \in T(T_{c(t_0)}M\slaz), 
\end{eqnarray*}
and since $t\mapsto J(t)$ is regular, we also have
\begin{eqnarray*}
  \kappa_2\circ J'(t_0) \in T(T_{c(t_0)}M\slaz).
\end{eqnarray*}
Since $c(t_0)\in \Sigma$, equation \eqref{pr:DFid_xx} implies that
$J'(t_0) = J'_1(t_0) = J'_2(t_0)$. Hence $J_1=J_2$. Since $J(0)\neq
0$, it follows that
\begin{eqnarray*}
  D\bigMap(\xi) = J_1'(0) 
                 = J_2'(0) 
                 = \kappa_2\circ D\bigMap\circ \kappa_2 ( \xi).
\end{eqnarray*}

\noindent\emph{Subcase B:} $J(t_0)$ arbitrary.\\
Let $j\colon \setR\times \varepsint \to TM$,
be the map
\begin{eqnarray*}
  j(t,s) &=& J(t) + sc'(t), \quad (t,s)\in \setR\times \varepsint.
\end{eqnarray*}
Now $j(\cdot,s)$ is a Jacobi field (with only isolated zeroes) for all
$s\in \varepsint$. If $s\in \varepsint\slaz$, Propositions
\ref{preJacobi_I} and \ref{preJacobi_II} imply that there exist Jacobi
fields $j_1(\cdot, s), j_2(\cdot, s) \colon \setR\to TM$ for
$\widetilde S$ such that\proofread{
Note that we do not assume that $j_i(\cdot,s)$ are smooth in $s$.}
\begin{eqnarray*}
  \partial_t j_1(t,s) &=& D\bigMap\circ \partial_t j(t,s), \quad  \quad\quad\quad\quad\mbox{when}\, t\in \setR \ \mbox{and} \  j(t,s)\neq 0, \\
  \partial_t j_2(t,s) &=& \kappa_2\circ D\bigMap\circ \kappa_2 \circ \partial_t j(t,s), \quad \mbox{when}\, t\in \setR.
\end{eqnarray*}
Let $\varepsilon>0$ be\proofread{
For $t_0$ we can find an $\varepsilon>0$ as follows:
\begin{enumerate}
\item If $J(t_0)=0$, then $j(t_0,s)=s c'(t_0)\neq 0$ for all $s\neq 0$, and
we can take any $\varepsilon>0$.
\item If $J(t_0)\neq 0$, then $j(t_0,0)=J(t_0)\neq 0$, and by continuity,
there exists a $\varepsilon>0$ such that $j(t_0,s)=J(t_0)+sc'(t_0)\neq 0$
for all $s\in \varepsint$. 
\end{enumerate}
} such that $j(t_0,s) \neq 0$ and $j(0,s) \neq 0$ for all $s\in
\varepsint\slaz$.  Then
\begin{eqnarray*}
  \partial_t j(t_0,s) \in T(T_{c(t_0)}M\slaz), \quad s\in \varepsint\slaz,
\end{eqnarray*}
and since $t\mapsto j(t,s)$ is regular, we also have
\begin{eqnarray*}
  \kappa_2\circ \partial_t j(t_0,s) \in T(T_{c(t_0)}M\slaz), \quad s\in \varepsint\slaz.
\end{eqnarray*}
Since $c(t_0)\in \Sigma$, equation \eqref{pr:DFid_xx} implies that
$\partial_t j(t_0,s) = \partial_t j_1(t_0,s) = \partial_t j_2(t_0,s)$
for all $s\in \varepsint\slaz$, so
$j_1(\cdot,s)=j_2(\cdot,s)$ for all $s\in \varepsint\slaz$. Let
$\Xi$ be the smooth curve $\Xi \colon \varepsint \to TTM\slaz$,
\begin{eqnarray*} 
  \Xi(s) &=& \partial_t j(t,s)|_{t=0}, \quad s\in \varepsint.
\end{eqnarray*} 
Then
\begin{eqnarray}
\label{bothSidesCont}
  D\bigMap\circ \Xi(s) &=& \kappa_2\circ D\bigMap\circ \kappa_2 \circ \Xi(s), \quad s\in \varepsint\slaz,
\end{eqnarray}
and equation \eqref{eq:toProve} follows 
since both sides of equation \eqref{bothSidesCont} are continuous
for $s\in \varepsint$ and since $\Xi(0)=\xi$.  
\commentOut{
Let $j\colon \setR\times \varepsint
\to TM$, be the map
\begin{eqnarray}
\label{Jdef}
  j(t,s) &=& J(t) + sc'(t), \quad (t,s)\in \setR\times \varepsint.
\end{eqnarray}
Now $t\mapsto j(t,s)$, $t\in \setR$, is a Jacobi field (with only
isolated zeroes) for all $s\in \varepsint$. 
Applying Propositions \ref{preJacobi_I} and
\ref{preJacobi_II} to each Jacobi field $j(\cdot,s)$, $s\in \varepsint$
and repeating the argument in Case 1b shows that
\begin{eqnarray}
\label{xxakljdas}
  D\bigMap\circ \Xi(s) &=& \wt \kappa_2\circ D\bigMap\circ \kappa_2 \circ \Xi(s), \quad s\in \varepsint\slaz,
\end{eqnarray}
where $\Xi$ is the smooth curve $\Xi \colon \varepsint \to TTM\slaz\cap
T(TM\slaz)$,
\begin{eqnarray*} 
  \Xi(s) &=& \partial_t j(t,s)|_{t=0}, \quad s\in \varepsint.
\end{eqnarray*} 
Equation \eqref{eq:toProve} follows since both sides in equation 
\eqref{xxakljdas} are continuous for $s\in \varepsint$ and 
since $\Xi(0)=\xi$. 
}
\end{proof}

\section{Descending maps and isometries}
\label{sec:GeoRie}
In this section we specialize Theorem \ref{mainThmXXSpray} to the case
when sprays $S$ and $\wt S$ are geodesic sprays of Riemann metrics. As
a result we obtain Theorem \ref{mainThmXX}, which gives sufficient
conditions for two Riemann metrics to be isometric.  It is not clear
whether Theorem \ref{mainThmXX} also hold for Finsler
metrics. However, the present proof uses that parallel transport is
norm-preserving for Riemann metrics. This result generalize to Berwald
metrics, but not to arbitrary Finsler metrics \cite[p.~89]{Shen2001b}.

The \emph{geodesic spray} of a (positive definite) Riemann metric $g$
is the spray with spray coefficients
\begin{eqnarray*}
  G^i(x,y) &=& \frac{1}{2} \Gamma^i_{ab}(x) y^a y^b,
\end{eqnarray*}
where $\Gamma^i_{jk}$ are the \emph{Christoffel symbols} associated
with $g$.

Suppose $c\colon I \to M$ is a geodesic for a Riemann metric $g$ and
$y\in T_{c(t)}M$ for some $t\in I$. Then there exists a unique curve
$V\colon I\to TM$ such that \emi $\pi_0\circ V=c$, \emii $V(t)=y$, and
\emiii $\nabla V=0$, where $\nabla V$ is covariant derivative induced
by $g$. We say that $V\colon I\to TM$ is the \emph{parallel transport}
of $y$ along $c$ and write $V(s) = P_{t\rightarrow s}(c)(y)$ for $s\in
I$. Thus $P_{t\rightarrow s}(c)$ is a linear map
$P_{t\rightarrow s}(c)\colon T_{c(t)}M \to T_{c(s)}M$.
If $\phi\colon M\to \widetilde M$ is a totally geodesic map
between Riemann manifolds and $c\colon I\to M$ is a geodesic, then
$\phi$ commutes with the parallel transport, so that
\cite{Vilms1970}
\begin{eqnarray}
\label{parCommute}
  (D\phi)(P_{t\rightarrow s}(c)(y)) &=& \widetilde P_{t\rightarrow s}( \phi\circ c)( D\phi(y)), \quad t,s\in I, \ y\in T_{c(t)}M.
\end{eqnarray}

\subsection{Isometric Riemann metrics}
\label{sec:Riemann}
Suppose   $\phi\colon M \to \widetilde M$ is a map and $g$ and
$\wt g$ are Riemann metrics on $M$ and $\wt M$, respectively. 
Then $\phi$ is an \emph{isometry} if for all $p\in M$, 
\begin{eqnarray}
  g(y,y) &=&\widetilde g(D\phi(y), D\phi(y)), \quad y\in T_pM.
\label{gisom}
\end{eqnarray}
Every isometry is a totally geodesic map
\cite[p.~232]{AbrahamMarsden:1994}. To prove Theorem
\ref{mainThmXX}, we will need the following converse result.
\begin{proposition} 
\label{guess22x}
Suppose $M$ and $\widetilde M$ are manifolds with Riemann metrics $g$
and $\wt g$, respectively. If $M$ is connected, $\phi$ is a
totally geodesic map $\phi\colon M\to \widetilde M$, and equation
\eqref{gisom} holds for one $p\in M$, then $\phi$ is an isometry.
\end{proposition}

\begin{proof}
For an open-closed argument, let
\begin{eqnarray*}
  A &=& \{ q\in M : g(y,y) = \widetilde g(D\phi(y), D\phi(y)) \,\, \mbox{for }\, y\in T_qM\}.
\end{eqnarray*}
By continuity, $A$ is closed, and by assumption, $A$ is non-empty. To
see that $A$ is open, let $q\in A$, and let $U\subset M$ be a normal
coordinate neighborhood around $q$.  If $r\in U$, then there exists a
geodesic $c\colon [0,1]\to M$ such that $c(0) = q$ and $c(1) = r$.
Then $\phi\circ c\colon [0,1]\to \wt M$ is also a geodesic.  Using
that parallel transport preserves Riemann norms, equation
\eqref{parCommute}, and that $q\in A$, it follows\proofread{ For $y\in
  T_rM$,
\begin{eqnarray*}
   g(y,y) &=& g( P_{1\to 0}(c)(y), P_{1\to 0}(c)(y) ) \\
               &=& \tilde g( (D\phi)\circ P_{1\to 0}(c)(y), (D\phi)\circ P_{1\to 0}(c)(y) ) \\
               &=& \tilde g( \tilde P_{1\to 0}(\phi\circ c)(D\phi(y)), \tilde P_{1\to 0}(\phi\circ c)(D\phi(y)) ) \\
               &=& \tilde g(D\phi( y), D\phi(y) ).
\end{eqnarray*}
} 
that $r\in A$. Thus $A$ is open, and $M=A$. 
\end{proof}
The next proposition shows that a Riemann metric is essentially
determined by its spray. This is a slight generalization of
Lemma 1 on page 242 in \cite{KobNom:1963}.
\begin{proposition} 
\label{SdetG}
Suppose $g$ and $\wt g$ are Riemann metrics on a connected manifold
$M$.  If $g$ and $\wt g$ have the same geodesic spray and $g=\wt g$ on
$T_pM$ for one $p\in M$, then $g=\wt g$.
\end{proposition}
\begin{proof}
  This follows by taking $M=\wt M$ and $\phi=\operatorname{id}$ in
  Proposition \ref{guess22x}. 
\end{proof}


\begin{theorem}
\label{mainThmXX}
Suppose $M$ is a smooth manifold $M$ with $\dim M \ge 2$ and $F$ is a
smooth map $\bigMap\colon TM\slaz\to TM\slaz$.  Furthermore, suppose
that $g$ and $\widetilde g$ are complete Riemann metrics on $M$ such
that $g$ has a trapping hypersurface $\Sigma\subset M$, and
\begin{enumerate}
\item $\bigMap$ maps Jacobi fields for $g$ into Jacobi fields for
  $\widetilde g$ (see Section \ref{sec:StrJacobi}),
\label{pr_g:jacobi}
\item for all $p\in \Sigma$, 
\begin{eqnarray*}
  S(y) &=& \wt S(y), \quad y\in T_pM\slaz,\\
  D\bigMap(\xi) &=& \xi,\quad\ \  \quad \xi\in T(T_pM\slaz),
\end{eqnarray*}
where $S$ and $\wt S$ are geodesic sprays induced by $g$ and $\wt g$,
respectively,
\item for one $p\in M$, 
\begin{eqnarray*}
  g(y,y) &=& \wt g(F(y), F(y)), \quad y\in T_pM\slaz.
\end{eqnarray*}
\end{enumerate}
Then there exists a smooth map $\phi\colon M\to M$ such that
\begin{enumerate}
\item $\bigMap=(D\phi)|_{TM\slaz}$, 
\item $\phi$ is an isometry (from $g$ to $\widetilde g$).
\end{enumerate} 
What is more, if $\bigMap$ is a diffeomorphism, then $\phi$ is a
diffeomorphism.
\end{theorem}

\begin{proof} 
  This follows from Theorem \ref{mainThmXXSpray} and Proposition
  \ref{guess22x}.
\end{proof}

\subsection*{Acknowledgements}
I.B. has been supported by 
   grant ID 398 from the Romanian Ministry of Education.
M.D. has been supported by 
   Academy of Finland Centre of Excellence programme 213476,
   the Institute of Mathematics at the Helsinki University of Technology,
and
   Tekes project MASIT03 --- Inverse Problems and Reliability of Models. 

\appendix
\section{Proof of Proposition \ref{preJacobi_I}}
\label{appZero}
For an outline of the proof below, see Section \ref{sec:StrJacobi} 

\begin{proof}[Proof of Proposition \ref{preJacobi_I}.]
  We can find a $t_0\in \setR$ such that $J$ restricts to a Jacobi
  field $J\colon I_0 \to TM\slaz$ without zeroes where $I_0\subset
  \setR$ is an neighborhood of $t_0$. Then $\wt J = \bigMap\circ J$
  defines a Jacobi field $\wt J\colon I_0\to T\wt M\slaz$ without
  zeroes. Since $\wt S$ is complete, $\wt S^c$ is complete
  \cite{BucataruDahl:2008}, and Jacobi field $\wt J$ extends into a
  Jacobi field $\wt J\colon \setR\to T\wt M$. For an open-closed
  argument, let $A=A_e \cup A_0$, where
\begin{eqnarray*}
  A_e &=& \{ t\in \setR : J(t)\neq 0 \ \ \mbox{and} \ \ \wt J'(t) = D\bigMap \circ J'(t)\},\\
  A_0 &=& \{ t\in \setR : J(t)= 0 \ \ \mbox{and} \ \ (t-\varepsilon, t)\cup (t,t+\varepsilon) \subset A_e \ \mbox{for some }\ \varepsilon>0 \}.
\end{eqnarray*}
Set $A$ is non-empty since $I_0 \subset A_e$.  
To see that $A$ is open, let us first note that $A_e$ is open since
$\bigMap$ maps Jacobi fields without zeroes to Jacobi fields without
zeroes and Jacobi fields are uniquely determined by their tangent at
one point. Also, if $t\in A_0$, then $t$ has a neighborhood $N\subset
\setR$ such that $N\setminus \{t\}\subset A_e$.

To see that $A$ is closed, let $t_i\in A$ be a sequence such that
$t_i\to \tau$ for some $\tau\in \setR$. Let us show that $\tau \in A$.
By Remark
\ref{zeroJacFieldRem}, we may assume\proofread{
  Since $A_e\cap A_0=\emptyset$, we may divide sequence $t_i$ into two
  (possibly finite) sub-sequences $\{e_i\}$ and $\{z_i\}$ such that
  $e_i\in A_e$ and $z_i\in A_0$.  There are two alternatives:
\begin{enumerate}
\item If $\#\{e_t\}=\infty$, we are done; sequence $e_i$ is a sequence
  in $A_e$ and $e_i\to \tau$.
\item If $\#\{e_t\}<\infty$, then $\#\{ z_i\}=\infty$ and we can
  divide sequence $z_i$ into two sub-sequences $\{\tau_i\}$ and
  $\{\zeta_i\}$ such that $\tau_i = \tau$ and $\zeta_i \neq
  \tau$. Again there are two alternatives:
\begin{enumerate}
\item[\emph{(a)}] $\#\{\zeta_i\}=\infty$ is not possibly by Remark
  \ref{zeroJacFieldRem}.
\item[\emph{(b)}] $\#\{\zeta_i\}<\infty$ implies that $\#\{\tau_i\}\ge
  1$ whence $\tau = t_i$ for some $i$ and $\tau \in A$.
\end{enumerate}
\end{enumerate}
}
that all $t_i\in A_e$.  
If $J(\tau)\neq 0$, then $\tau \in A_e$ by\proofread{
Since $t_i\to \tau$, $\wt J'(\tau)= \lim \wt J'(t_i) = \lim DF\circ J'(t_i)=DF\circ J'(\tau)$.
}
continuity.
If $J(\tau)=0$, we show that $\tau \in A_0$.  This is
straightforward\proofread{
  By Remark \ref{zeroJacFieldRem}, there exists an $\varepsilon>0$
  such that $J$ is non-zero on $(\tau-\varepsilon,
  \tau+\varepsilon)\setminus\{\tau\}$. By assumption there are
  $t_\pm\in A_e$ such that $t_-\in (\tau-\varepsilon, \tau)$ and
  $t_+\in (\tau, \tau+\varepsilon)$.  By uniqueness, we then have
  $F\circ J=K$ on $(\tau-\varepsilon, \tau)$ and on $(\tau,
  \tau+\varepsilon)$, and since $J(\tau)=0$, we have $\tau \in A_0$.
}
to check using uniqueness if an arbitrary neighborhood of $t$ contains
$t_i$:s on both sides of $\tau$. Let us assume that $t_i<\tau$ for all
$i\ge 1$. (The case $t_i>\tau$ is analogous.)

Let $j\colon I \times \varepsint\to TM$ be the map obtained by
applying Lemma \ref{lemma:puncturedJacobiVariation} below to $J$. Then
$\tau \in I$ and $j(t,s)\neq 0$ on for $(t,s)\neq (\tau,0)$. Let $\wt
j$ be the map
\begin{eqnarray*}
\wt j\colon (I\times \varepsint )\setminus\{(\tau,0)\} &\to& T\wt M\slaz \\
  \wt j(t,s) &=& \bigMap\circ j(t,s).
\end{eqnarray*}
For each $s\in \varepsint\slaz$, $\wt j(\cdot,s)\colon I\to T\wt
M\slaz$ is a Jacobi field without zeroes, and for $s=0$, $\wt
j(\cdot,0)\colon I_\pm\to T\wt M\slaz$ are Jacobi fields without zeroes,
where
$$
  I_+ \ = \ \{t \in I : t>\tau\},
  \quad 
  I_- \ = \ \{t \in I : t<\tau\}.
$$
We know that $\wt J=\wt j(\cdot,0)$ on $I_-$, and $\tau\in A_0$ follows if
$\wt J=\wt j(\cdot,0)$ on $I_+$. If $\wt \Phi^c_t$ is the flow of $\wt S^c$, and 
$t_-\in I_-$, then for $t_+\in I_+$ we have
\begin{eqnarray*}
  \wt j(t_+,0) &=& \lim_{s\to 0} \wt j(t_+,s) \\
         &=& \lim_{s\to 0} \wt \pi_1\circ \wt \Phi^c_{t_+-t_-}( \partial_t \wt j(t_-,s)) \\
         &=& \wt \pi_1\circ \wt \Phi^c_{t_+-t_-}( \wt J'(t_-)) \\
         &=& \wt J(t_+). \qedhere
\end{eqnarray*}
\end{proof}

\begin{lemma} 
\label{lemma:puncturedJacobiVariation}
Suppose $\dim M\ge 2$, $J\colon \setR \to TM$ is a Jacobi field for
spray $S$, and $\tau\in \setR$ is an isolated zero for $J$. Then
$\tau$ has a neighborhood $I\subset \setR$, and there exists a map
$j\colon I\times \varepsint \to TM$ such that
\begin{enumerate}
\item $j(t,0)=J(t)$ for $t\in I$,
\item $t\mapsto j(t,s)$, $t\in I$, is a Jacobi field for all $s\in \varepsint$, 
\item $j(t,s)\neq 0$ if $(t,s)\neq (\tau,0)$.
\end{enumerate}
\end{lemma}

\begin{proof} We may assume that $\tau = 0$. Let $c$ be geodesic
  $c\colon \setR\to M$, $c=\pi_0\circ J$.  In local coordinates, we
  have $J'(0)=(x(0), 0, \dot{x}(0), \dot{J}(0))$, and let $\xi\in
  T_{c(0)}M\slaz$ be vector $\xi=(x(0),\dot J (0))$. Then there exists
  an auxiliary Riemann metric on $M$ such that $g(\xi,\xi)=1$, and
  since $\dim M\ge 2$, there exists a non-zero vector $v\in
  T_{c(0)}M\slaz$ such that $g(v,\xi)=0$ and $g(v,v)=1$.  Let $K$ be a
  Jacobi field $K\colon I\to TM$ determined by $K'(0) = (x(0),v,\dot
  x(0), 0)$, and let $j$ be the map $j\colon I\times \setR \to TM$
  defined as
\begin{eqnarray*}
  j(t,s) &=& J(t) + s K(t), \quad (t,s)\in I\times \setR.
\end{eqnarray*}
Now \emi and \emii are clear. For \emph{(iii)}, let us shrink $I$ such
that $c\colon I\to M$ is contained in the domain of coordinates
$x^i$. Then $j$ has local expression $j(t,s)=(x^i(t), j^i(t,s))$, and
\begin{eqnarray*}
  j^i(t,s) = \xi^i t + v^i s + R^i(t,s),  \quad (t,s)\in I\times \setR,
\end{eqnarray*}
with remainder terms $R^i(t,s)={o}(\sqrt{t^2+s^2})$.
For curve $j_0\colon I\times \setR \to T_{c(0)}M$,
\begin{eqnarray*}
  j_0(t,s) &=& (x^i(0), j^i(t,s)),
\end{eqnarray*}
the Cauchy--Schwarz inequality yields
\begin{eqnarray*}
  g(j_0(t,s),j_0(t,s))&=& t^2 + s^2 + o(t^2 + s^2), \quad
  (t,s) \in I\times \setR,
\end{eqnarray*}
and 
\proofread{ 
For a contradiction, suppose that $(t_i, s_i)$ is a sequence such that 
\begin{enumerate}
\item $(t_i, s_i)\to (0,0)$, 
\item 
$g(j_0(t_i,s_i),j_0(t_i,s_i))< \frac{1}{2}(t_i^2 + s_i^2)$.
\end{enumerate}
Then there exists a function $R(t,s)=o(t^2 + s^2)$ such that $t_i^2 +
s_i^2< R(t_i,s_i)$ for all $i\ge 1$. However, this is a contradiction
since $\lim_{i\to \infty} \frac{\operatorname{LHS}}{t_i^2 + s_i^2} =
1$, but $\lim_{i\to \infty} \frac{\operatorname{RHS}}{t_i^2 + s_i^2} =
0$.  } 
\emiii follows since we can find $\varepsilon>0$ such that
\begin{eqnarray*}
  g(j_0(t,s),j_0(t,s))&\ge& \frac{1}{2}(t^2 + s^2), \quad (t,s)\in \varepsint^2.\qedhere
\end{eqnarray*}
\end{proof}

\providecommand{\bysame}{\leavevmode\hbox to3em{\hrulefill}\thinspace}
\providecommand{\MR}{\relax\ifhmode\unskip\space\fi MR }
\providecommand{\MRhref}[2]{%
  \href{http://www.ams.org/mathscinet-getitem?mr=#1}{#2}
}
\providecommand{\href}[2]{#2}


\end{document}